\documentclass[12pt]{amsart}

\newtheorem{theorem}{Theorem}[section]
\newtheorem{proposition}[theorem]{Proposition}
\newtheorem{lemma}[theorem]{Lemma}
\newtheorem{cor}[theorem]{Corollary}

\theoremstyle{definition}

\newtheorem*{rem}{[Remark]}

\newtheorem{exam}[theorem]{Example}

\newcommand{\nn}{\ensuremath{[n]}}
\newcommand{\nset}{\ensuremath{{{[n]}\choose {2}}}}
\newcommand{\ij}{\ensuremath{\{i,j\}}}
\newcommand{\phai}{\ensuremath{\varphi_{g,k}}}
\newcommand{\dele}[1]{\ensuremath{\Delta^e({#1})}}

\newcommand{\dels}[1]{\ensuremath{\Delta^s({#1})}}
\newcommand{\pusai}{\ensuremath{\phi_{g,k}}}
\newcommand{\gin}{\ensuremath{\mathrm{gin}}}
\newcommand{\Gin}{\ensuremath{\mathrm{Gin}}}
\newcommand{\GIN}{\ensuremath{\mathrm{GIN}}}

\newcommand{\delc}[1]{\ensuremath{\Delta^c({#1})}}
\newcommand{\shift}{\ensuremath{\mathrm{Shift}}}
\newcommand{\GL}{\ensuremath{GL_n (K)}}
\newcommand{\init}{\ensuremath{\mathrm{in}}}
\newcommand{\tuua}{\ensuremath{\tau_{g,k}^A}}
\newcommand{\tuub}{\ensuremath{\tau_{g,k}^B}}

\newcommand{\shif}{\ensuremath{\mathrm{Shift}_{ij}(\ensuremath{G})}}

\def\cocoa{{\hbox{\rm C\kern-.13em o\kern-.07em C\kern-.13em o\kern-.15em A}}}

\begin{document}
\title{Algebraic shifting of finite graphs}
\author{Satoshi Murai}
\address{
Department of Pure and Applied Mathematics\\
Graduate School of Information Science and Technology\\
Osaka University\\
Toyonaka, Osaka, 560-0043, Japan\\
}
\email{s-murai@ist.osaka-u.ac.jp
}
\maketitle

\begin{abstract}
%Algebraic shifting is a correspondence which associate to a graph $G$ an another graph $\Delta(G)$ with special conditions.
%There are two main variants exterior shifting and symmetric shifting.
In the present paper,
exterior algebraic shifting and symmetric algebraic shifting of bipartite graphs
and chordal graphs are studied.
First, we will determine the symmetric algebraic shifted graph of complete bipartite graphs.
It turns out that, for $a>3$ and $b>3$,
the exterior algebraic shifted graph of the complete bipartite graph $K_{a,b}$ of size $a,b$ is different from
the symmetric algebraic shifted graph of $K_{a,b}$.
Second, we will show that the exterior algebraic shifted graph of any chordal graph $G$
is coincident with the symmetric algebraic shifted graph of $G$.
\end{abstract}

\section*{Introduction}
Let $G=(\nn,E(G))$ be a finite graph on the vertex set $\nn=\{1,2,\dots,n\}$ with the edge set $E(G)$.
Throughout this paper, we assume that all graphs have no loops and no multiple edges.
A graph $G$ is called \textit{shifted} if, for any edge $\{i,j\}\in E(G)$ and
for any integers $i'\leq i$ and $j'\leq j$, one has $\{i',j'\}\in E(G)$.
The shifting operation of graphs is an operation which associate with a graph $G$ a shifted graph $\Delta (G)$.
Shifting operations are first considered by Erd\"os, Ko and  Rado \cite{EKR}.
Their shifting operation $G \to \delc G$ is called combinatorial shifting.
On the other hand, algebraic shifting was introduced by Kalai.
Main variants of algebraic shifting are exterior algebraic shifting $G \to \dele G$
and symmetric algebraic shifting $G \to \dels G$.
The definition of algebraic shifting will be given in \S 1.
Shifting operations are defined for simplicial complexes, however,
we only consider finite graphs in the present paper.

Algebraic shifting of graphs was first studied in \cite{KA-h}.
The connectivity of $\dele G$ was mainly studied in it.
In the present paper, we are interested in the comparison between exterior algebraic shifting and symmetric algebraic shifting.
Our main problem is to analyze when $\dele G$ and $\dels G$ are equal and when they are different.
The complete bipartite graph $K_{3,3}$ of size 3,3 is an example
whose exterior algebraic shifted graph $\dele {K_{3,3}}$ and symmetric algebraic shifted graph $\dels {K_{3,3}}$
are different (this example was found by Kalai).
Also, we can easily show that  $\dele G= \dels G$ if $G$ is tree.
In the present paper, to extend these two examples,
we will study algebraic shifting of bipartite graphs and chordal graphs.

First, we consider  bipartite graphs in \S 2.
Let $K_{a,b}$ be the complete bipartite graph of size $a,b$.
Our result determines the symmetric algebraic shifted graph of complete bipartite graphs.
Kalai \cite{KA-h} determined the exterior algebraic shifted graph of complete bipartite graphs.
Kalai's result together with our result explains the difference between $\dele {K_{a,b}}$ and $\dels {K_{a,b}}$.
In fact, 
we have  $\dele {K_{a,b}} \ne \dels {K_{a,b}}$ if $a\geq 3$ and $b \geq 3$.
Also, we will show that if $G$ is a bipartite graph and has sufficiently many edges, then one has $\dele G \ne \dels G$.

Second, we consider chordal graphs.
A finite graph $G$ is called \textit{chordal} if every cycle in $G$ of length $> 3$ has a chord.
We will show that every chordal graph $G$ satisfies $\dele G=\dels G$.
This fact will be given in \S 3.

At last, we will consider combinatorial shifting of chordal graphs in \S 4.
In general, computations of algebraic shifting are rather difficult.
No effective method to compute algebraic shifting is known.
%In fact, it is really hard to compute algebraic shifting.
On the other hand, combinatorial shifting is easily computed.
However, combinatorial shifting does not behave nicely.
For example, combinatorial shifting is not even uniquely determined.
We are interested in when we can realize  $\dele G$ as a combinatorial shifted graph.
In general,
there exist a simplicial complex $\Gamma$ such that $\delc \Gamma \ne \dele \Gamma$ for an arbitrary combinatorial shifted complex
$\delc \Gamma$ of $\Gamma$.
% for a simplicial complex $\Gamma$,
%there does not always exists a combinatorial shifted complex $\delc \Gamma $ such that $\delc \Gamma = \dele \Gamma$
(such examples are introduced in \cite[\S 6.2]{KA} and in \cite[{\S 4}]{H-M}). 
However, we do not have an example of a finite graph $G$ which satisfies  $\delc G \ne \dele G$
for an arbitrary combinatorial shifted graph $\delc G$.
We will show that
there is a combinatorial shifted graph $\delc G$ such that $\delc G = \dele G$ for every chordal graph $G$.
Furthermore,
we will give an algorithm to compute the combinatorial shifted graph $\delc G$ with $\delc G=\dele G$
when $G$ is chordal.
Also, all chordal graphs, whose combinatorial shifted graph is unique,
will be determined in Proposition \ref{chordunique}.

\section{algebraic shifting}

Let $K$ be a field of characteristic $0$ and
$R =K[x_1,x_2,\dots,x_n]$
the polynomial ring in $n$ variables over a field $K$ with each $\deg(x_i)=1$.
Fix a term order $\prec$ on $R$.
For an ideal $I\subset R$,
let $\init_\prec (I)$ denote the initial ideal of $I$ w.r.t. the term order $\prec$.

Let $\GL$ be the general linear group with coefficients in $K$.
For any $g=(a_{ij}) \in \GL$,
define the ring automorphism $g:R \to R$ by 
$$ g (f(x_1,x_2,\dots,x_n))=f(\sum_{i=1}^n a_{i1}x_i,\sum_{i=1}^n a_{i2}x_i,\dots,\sum_{i=1}^n a_{in}x_i)
\ \ \mbox{for all } f\in R.$$
For an ideal $I$,
we write $g(I)=\{g(f):f\in I\}$.

\begin{theorem}[Galligo, Bayer and Stillman] \label{gin}
Fix a term order $\prec$.
% satisfying $x_1 \prec x_2 \prec \dots \prec x_n$.
For each homogeneous ideal $I \subset R$,
there is a nonempty Zariski open subset $U\subset \GL$ such that
$\init_{\prec} ( g (I))$ is constant for all $g \in U$.
%Furthermore $\init_{<} (\varphi (I))$ with $\varphi \in U$ is strongly stable.
\end{theorem}

The above monomial ideal $\init_{\prec}(g (I))$ with $g \in U$ is called the
{\em generic initial ideal of $I$ with respect to the term order $\prec$},
and will be denoted $\gin_{\prec} (I)$.
Let $\prec_{rev}$ be the degree reverse lexicographic order induced by $x_1\prec x_2\prec \cdots \prec x_n$.
In other words,
for monomials $u=x_1^{a_1} \cdots x_n^{a_n}$ and $v = x_1^{b_1} \cdots x_n^{b_n}$,
define $u \prec_{rev} v$ if $\deg(u)<\deg(v)$ or $\deg(u)=\deg(v)$ and 
the leftmost nonzero entry of $(b_1-a_1,\dots,b_n -a_n)$ is negative.
In the present paper,
we
% only consider generic initial ideals w.r.t. the degree reverse lexicographic order and
write $\gin (I)= \gin_{\prec_{rev}} (I).$

We say that a property $\mathcal{P}$ holds for a generic matrix $g \in \GL$ if there is a nonempty Zariski open subset $U\subset \GL$
such that $\mathcal{P}$ holds for all $g \in U$.

\medskip
\noindent\textbf{(Symmetric algebraic shifting)}
Let $[n]=\{1,2,\dots,n\}$ and
${{[n]}\choose {2}}=\{ \ij \subset [n] : \ i\ne j \}$.
Let $G=([n],E(G))$ be a graph on $\nn$ with the edge set $E(G)\subset {{[n]}\choose {2}}$.
The \textit{Stanley--Reisner ideal} $I_G \subset R$ of $G$ is the ideal
generated by all squarefree monomials $x_ix_j$ with $\{i,j\}\not\in E(G)$
%$\{x_ix_j|\{i,j\}\not \in E(G)\}$
and all squarefree monomials of degree $3$.
Define
$$\GIN({G})=\{x_i x_j \in R:\ x_ix_j \not \in \gin(I_G)\}.$$
The \textit{symmetric algebraic shifted graph}  $\dels G$ of ${G}$
is the graph on $[n]$ defined by
$$E(\dels {G})=\{ \{i-1,j\} \subset [n] :x_i x_j \in \GIN({G}), 2 \leq i \leq j \leq n\}.$$
The operation ${G}\to \dels {G}$
is called \textit{symmetric algebraic shifting}.

\medskip
\noindent\textbf{(Exterior algebraic shifting)}
Let $V$ be a $K$-vector space of dimension $n$ with basis $e_1,\dots,e_n$
and $E=\bigoplus_{d=0}^n \bigwedge^d V$ the exterior algebra of $V$.
For a subset $S = \{ s_1, s_2,\dots,s_k \} \subset [n]$ with $s_1< s_2< \cdots < s_k$,
we write $e_S=e_{s_1}\wedge e_{s_2} \wedge \cdots \wedge e_{s_k}\in E$. 
The element $e_S$ is called the \textit{monomial} of $E$ of degree $k$.
In the exterior algebra,
the generic initial ideal $\Gin(J)$ of a graded ideal $J\subset E$
is defined similarly as in the case of the polynomial ring 
(\cite[Theorem 1.6]{AHH}).

Let $G$ be a graph on $[n]$.
The \textit{exterior face ideal} $J_G \subset E$ of $G$ is the ideal
generated by all monomials $e_i \wedge e_j$ with $\{i,j\} \not \in E(G)$
% $\{e_i\wedge e_j:\{i,j\}\not \in E(G)\}$
and all monomials in $E$ of degree $3$.
The \textit{exterior algebraic shifted graph}  $\dele G$ of ${G}$
is the graph on $[n]$ defined by 
$$E(\dele {G})=\{ \ij \subset [n]:e_i \wedge e_j \not\in \Gin(J_G)\}.$$
The operation ${G}\to \dele {G}$
is called \textit{exterior algebraic shifting}.

We recall some basic properties of generic initial ideals and algebraic shifting.

\begin{lemma} \label{tuketasi}
Let $I=\oplus_{d \geq 0} I_d$ be a homogeneous ideal,
where $I_d$ is the $d$-th homogeneous component of $I$.
Then
\begin{itemize}
\item[(i)]
$I$ and $\gin (I)$ have the same Hilbert function.
That is, one has $\dim_K (I_d)= \dim_K (\gin(I)_d)$ for all $d \geq 0$.
\item [(ii)]
If $I\subset J$ are homogeneous ideals,
 then one has $\gin (I) \subset \gin(J)$.  
\end{itemize}
\end{lemma}

\begin{lemma}[{{\cite[\S 8]{Her} or \cite[\S 2]{KA}}}] \label{itumono}
Let $G$ and $G'$ be graphs on $[n]$. Then
\begin{itemize}
\item[(i)] $\dele G$ and $\dels G$ are shifted;
\item[(ii)] $|E(G)|=|E(\dele G)| =|E( \dels G)|$;
\item[(iii)] If $E(G) \subset E(G')$ then  $E(\dele G) \subset E( \dele {G'})$ and $E(\dels G) \subset E(\dels {G'})$.
\end{itemize}
\end{lemma}

We recall the technique to compute exterior algebraic shifting developed in \cite{KA-h}.
Let  $g=(a_{ij}) \in GL_n(K)$ and  $1 \leq k \leq n$.
Define a map $\phai:\nset \to \bigoplus_{i=1}^k V$ by
\[
\phai (i,j)=(a_{1j}e_i-a_{1i}e_j,\dots,a_{kj}e_i-a_{ki}e_j) \in \bigoplus_{i=1}^k V,
\]
where $1 \leq i <j \leq n$.
For a graph $G$ on $[n]$,
let $\phai(G)=\mathrm{span}\{ \phai (i,j): \{i,j\} \in E(G)\}$
be the vector space spanned by all vectors $\phai(i,j)$ with $\{i,j\} \in E(G)$.
Set $r_{g,k}(G) = \dim_K (\phai(G))$.

For a shifted graph $G$ on $[n]$ and for an integer $1 \leq k \leq n-1$, define 
$$m_{ k}(G)=|\{\ij \in E(G): \min\ij = k\}|$$
 and 
$$m_{\leq k}(G)=|\{\ij \in E(G): \min\ij\leq k\}|,$$
where $|A|$ denote the cardinality of a finite set $A$.
Note that $m_1(G)> m_2 (G) >\cdots >m_{n-1}(G)$ since $G$ is shifted.
Then $m_{\leq 1}(G),\dots,m_{\leq n-1}(G)$ determine a shifted graph $G$ by the relations
\begin{eqnarray}
E(G)=\{ \{i,j\} \subset [n]: 1 \leq i \leq n-1,\ i<j \leq i+m_i(G)\} \label{tondemo}
\end{eqnarray}
and 
$$m_k(G) = m_{\leq k}(G) -m_{\leq k-1}(G),$$
where we let $m_0(G)=0$.
%Thus knowing $\dele G$ is equivalent to knowing $m_{ \leq 1}(\dele G),\dots, m_{\leq n}( \dele G)$.

\begin{lemma}[{ {\cite[Lemma 7.1]{KA-h}}}] \label{exdet}
Let $G$ be a graph on $[n]$.
For a generic matrix $g \in \GL$, one has
$$m_{\leq k } (\dele G) = r_{g,k}(G)\ \ \ \mbox{ for }\ k=1,2,\dots,n-1.$$
\end{lemma}

We also recall similar technique for symmetric algebraic shifting appeared in \cite[\S 8]{KA-h}.
Let  $g=(a_{ij}) \in GL_n(K)$ and  $1 \leq k \leq n$.
For integers $1 \leq i \leq j \leq n$, define
$$\pusai (i,j) = (a_{1j}e_i+a_{1i}e_j,\dots,a_{kj}e_i+a_{ki}e_j)\in \bigoplus_{i=1}^k V.$$
For a finite graph $G$ on $[n]$,
let 
\begin{eqnarray*}
\pusai(G)&=&\mathrm{span}\{ \pusai (i,j): \{i,j\} \in E(G)\}, \\
\pusai([n])&=&\mathrm{span}\{ \pusai (i,i): 1 \leq i \leq n\}, \\
\pusai(G^+)&=&\mathrm{span}\{ \pusai (i,j): \{i,j\} \in E(G) \mbox{ or } 1 \leq i = j \leq n \}
\end{eqnarray*}
and $s_{g,k}(G) = \dim_K (\pusai(G^+))$.

To give an analogue of Lemma \ref{exdet},
we recall some fundamental facts.

\begin{lemma} \label{initial}
Let $\mathcal{M}$ be the set of monomials in $R$,
$\prec$ a term order and $I \subset R$ an ideal.
Let $ u \in R$ be a monomial and $\tilde u$ its image in $R/I$.
Then
$u \not\in \init_{\prec}(I) $ if and only if
$\tilde u \not\in \mathrm{span}\{ \tilde v: v \prec u,\ v \in \mathcal{M}\}$.
\end{lemma}

The above lemma is well known and easy.
Thus we omit the proof.

\begin{lemma} \label{tyoutenn}
Let $n \geq m >0$ be positive integers.
Let $G=([n],E(G))$ and $G'=([m],E(G'))$ be graphs with $E(G)=E(G')$.
Then, one has $E(\dels G) =E(\dels {G'})$.
\end{lemma}

\begin{proof}
We may assume $n=m+1$.
Let $R=K[x_1,\dots,x_n]$.
Note that $I_{G'}$ is the ideal of $K[x_1,\dots,x_m]$.
Consider the ideal $J=I_{G'}R+x_nR$ of $R$.
Since $\gin(x_n R)= x_nR$, Lemma \ref{tuketasi} says that $x_n \in J$.
We claim $\gin(J)=\gin(I_{G'})R + x_nR$. 

For $\varphi=(a_{ij}) \in \GL$ with $a_{nn}\ne 0$,
let $\tilde \varphi=(a_{ij}-\frac {a_{in}} {a_{nn}} a_{nj})_{1 \leq i,j \leq m} \in \mathrm{GL}_m(K)$.
If $U \subset \GL$ is a nonempty Zariski open subset,
then $\{ \tilde \varphi: \varphi=(a_{ij}) \in U, \ a_{nn}\ne 0\}$ contains a nonempty Zariski open subset of $\mathrm{GL}_m (K)$.
Thus there exists $\varphi \in \GL$ such that
$\gin(J)= \init_{\prec_{rev}}(\varphi(J))$ and $\gin (I_{G'}) = \init_{\prec_{rev}}(\tilde \varphi(I_{G'}))$.

On the other hand,
for $i = 1,2,\dots,m$,
we have $\tilde \varphi(x_i)= \varphi(x_i-\frac {a_{in} } {a_{nn}}x_n)$.
Let $x_1^{a_1} \cdots x_m^{a_m} \in I_{G'}$.
Then we have
$$\tilde \varphi (x_1^{a_1} \cdots x_m^{a_m})
=\varphi(x_1-\frac {a_{1n}}{a_{nn}} x_n)^{a_1} \cdots \varphi(x_m- \frac{a_{mn}} {a_{nn}} x_n)^{a_m} 
\in \varphi(x_1^{a_1} \cdots x_n^{a_n})R +\varphi(x_n)R.$$
Thus we have $\tilde \varphi(x_1^{a_1} \cdots x_n^{a_n}) \in \varphi(J)$ for all $x_1^{a_1} \cdots x_n^{a_n} \in I_{G'}$.
Thus $\gin(I_{G'}) \subset \gin(J)$ and
$\gin(I_{G'})R + x_nR \subset \gin(J)$ since $x_n \in \gin (J)$.
Since $\gin(I_{G'})R + x_nR$ and $\gin (J)$ have the same Hilbert function,
we have $\gin(I_{G'})R +x_nR=\gin (J)$.

Now, since $I_G \subset I_{G'}R +x_nR$, we have $\gin (I_G) \subset \gin (I_{G'})R + x_nR$
and $\GIN(G) \supset \GIN(G')$.
Thus we have $E(\dels G) \supset E(\dels {G'})$.
Then, by Lemma \ref{itumono} (ii), we have $E(\dels G)= E(\dels {G'})$.
\end{proof}

Let $G$ be a graph on $[n]$.
For an integer $1 \leq k \leq n$,
write 
$$m_{\leq k} (\GIN(G)) = |\{ x_i x_j \in \GIN(G): \min\{ i,j\} \leq k\}|.$$

\begin{lemma}\label{symmetric}
Let $G$ be a graph on $[n]$.
Then
\begin{itemize}
\item[(i)] $m_{\leq 1}(\GIN(G))=n$;
\item[(ii)] $m_{\leq k}(\dels G) = m_{\leq k+1}(\GIN(G) )-n$\ \  for all $k=1,2,\dots,n-1$.
\end{itemize}
\end{lemma}

\begin{proof}
(i) By Lemma \ref{itumono},
we have 
$|\{x_ix_j\in \GIN(G): 2 \leq i \leq j \leq n\}|=|E(G)|.$
On the other hand,
since $I_G$ and $\gin(I_G)$ have the same Hilbert function,
we have
$$|\GIN(G)|=|\{x_ix_j\in \mathcal{M}:x_ix_j \not \in I_G\}|=|E(G)|+n.$$
Thus we have $m_{\leq 1} (\GIN(G))=n$.

(ii) By the definition, we have
$m_{\leq k}(\dels G) =m_{\leq k+1}(\GIN(G))-m_{\leq 1}(\GIN(G)).$
Then the assertion immediately follows from (i). 
\end{proof}

\begin{lemma} \label{GIN}
Let $G$ be a graph on $[n]$.
For a generic matrix $g\in \GL$,
one has
$$s_{g,k}(G)=m_{\leq k}(\GIN(G))\ \ \ \mbox{ for } k=2,3,\dots,n.$$
In particular, one has $m_{\leq k }( \dels G) =s_{g,k+1}(G)-n$ for $k=1,2,\dots,n-1$.
\end{lemma}

\begin{proof}
Let $g \in \GL$ be a matrix with $\gin(I_G)= \init_{\prec_{rev}} (g(I_G))$.
Set $g^{-1}= (a_{ij})$ and $y_i = g(x_i)$ for $i =1,\dots,n$.
Then $x_i=g^{-1}(y_i) = \sum_{l=1}^n a_{li}y_l$ for $i= 1,\dots,n$
and
\begin{eqnarray}
\quad \quad x_px_q = \sum_{1 \leq i < j \leq n} (a_{ip}a_{jq}+a_{iq}a_{jp})y_iy_j + \sum_{i=1}^n a_{ip}a_{iq} y_i^2\ \ \ \ \ \ 
\mbox{ for all }1 \leq p,q \leq n. \label{UFO}
\end{eqnarray}
For integers $1 \leq i \leq j \leq n$ and $1 \leq p \leq q \leq n$,
let
\begin{eqnarray*}
\alpha_{ (p,q)(i,j)}=
\left\{
\begin{array}{l}
a_{ip}a_{jq} + a_{iq}a_{jp}, \ \ \ \ \mbox{ if } i \ne j,\\
a_{ip}a_{iq}, \ \ \ \ \hspace{43pt}\mbox{ if } i=j.
\end{array}
\right.
\end{eqnarray*}
Consider the ${n+1 \choose 2} \times \{|E(G)|+n\}$ matrix $M=(\alpha_{(p,q)(i,j)})$
%_{ \ij \in E(G) \mbox{ or }i=j,\ 1 \leq p \leq q \leq n}$
whose lows are indexed by $(p,q)$ with $1\leq p\leq q \leq n$
and whose columns are indexed by $(i,j)$ with $\{i,j\} \in E(G)$, where $i<j$, or $1 \leq i=j \leq n$.
For $1 \leq k \leq n$,
let $M_k$ be the submatrix of $M$ whose lows are indexed by $(p,q)$ with $1 \leq p \leq k$.
%and with $1 \leq p \leq q \leq n$. 

Let $x_p x_q$ be a monomial with $1 \leq p \leq q \leq n$.
We write $\widetilde{x_px_q}$ for the image of $x_px_q$ in $R/ g(I_G)$. 
Since $g(I_G)$ is generated by the elements $y_i y_j$ with $\{i,j\} \not \in E(G)$ together with
all squarefree monomials in $K[y_1,\dots,y_n]$ of degree $3$.
The form (\ref{UFO}) says that
$\widetilde{x_px_q}$ is coincident with the $(p,q)$-th low vector of $M$.

The definition of $\prec_{rev}$ says that 
$x_px_q \preceq_{rev} x_kx_n$ if and only if $1 \leq p \leq k$.
Since $\gin(I_G)= \init_{\prec_{rev}}(g(I_G))$,
Lemma \ref{initial}
together with the definition of $\prec_{rev}$ 
says that
$m_{\leq k}(\GIN(G))$ is equal to the number of $K$-linearly independent elements in
$\{ \widetilde{x_px_q} \in R/ g(I_G): 1 \leq p \leq k \mbox{ and }  p \leq q \leq n\}.$
Since each $\widetilde{x_px_q}$ is coincident with $(p,q)$-th low vector of $M_k$,
we have
\begin{eqnarray}
m_{\leq k} (\GIN(G)) = \mathrm{rank}(M_k) \ \ \ \mbox{ for } k= 1,2,\dots,n. \label{sono4}
\end{eqnarray}

Let $W$ be a ${n+1 \choose 2}$-dimensional $K$-vector space whose basis are $e_{pq}$ with $1\leq p \leq q \leq n$
and $W_l$ the subspace of $W$ spanned by $e_{ll},e_{ll+1},\dots,e_{ln}$.
For integers $1 \leq i \leq j \leq n$ with $\ij\in E(G)$ or $i=j$,
the vector
$$f_k(i,j)=\sum_{{1 \leq p \leq k}\atop {p \leq q \leq n}}\alpha_{(p,q)(i,j)}e_{pq} \in \bigoplus_{l=1}^k W_l$$
is coincident with the column vector of $M_k$.
Thus we have
\begin{eqnarray}
\dim_K( \mathrm{span}\{f_k(i,j):\{i,j\} \in E(G) \mbox{ or } 1\leq i=j\leq n \} )=\mathrm{rank}(M_k). \label{ppc2}
\end{eqnarray}
For integers $1 \leq l \leq n$ and $1 \leq i \leq n$,
let 
$$v_{li}= \frac{1}{2} a_{il}e_{ll} + \sum_{t=l+1}^n (a_{it}e_{lt}) \in W_l.$$
Then $f_k(i,j)$ can be written of the form
$$f_k(i,j)=\sum_{l=1}^k (a_{il}v_{lj}+a_{jl}v_{li}).$$

Let $m= \max\{ \max \ij : \ij \in E(G)\}$.
When we consider $m_{ \leq k}(\GIN (G))$ with $k \geq 2$,
Lemma \ref{tyoutenn} says that we may assume $n$ is sufficiently larger than $m$.
Thus, for a generic matrix $g=(a_{ij}) \in \GL$ and for $1 \leq l \leq k$,
we may regard $v_{l1},v_{l2},\dots,v_{lm}$ as $K$-linearly independent vectors of $W_l$
which do not depend on $a_{st}\in K$ with $1 \leq s \leq m$.
Then, for a generic matrix $g\in \GL$ and for integers $1 \leq i,j \leq n$ and $2 \leq k \leq m$,
the vector $f_k(i,j)$ is coincident with $\phi_{g,k} (i,j)$.
%Thus, for $2 \leq k \leq m$, we have
%\begin{eqnarray*}
%&&\dim_K (\mathrm{span}\{(f_k(i,j)):\{i,j\} \in E(G) \mbox{ or } 1 \leq i=j \leq m\})\\
%&&=\dim_K (\mathrm{span}\{\pusai(i,j):\{i,j\} \in E(G) \mbox{ or } 1 \leq i=j \leq m\}) \\
%\end{eqnarray*}

On the other hand,
since $v_{1m+1},\dots,v_{1n}$ appear only in $f_{k}(m+1,m+1),\dots,f_k(n,n)$,
we have
\begin{eqnarray*}
&&\dim_K( \mathrm{span}\{(f_k(i,j)):\{i,j\} \in E(G) \mbox{ or } 1 \leq i=j \leq n\})\\
&&=\dim_K (\mathrm{span}\{(f_k(i,j)):\{i,j\} \in E(G) \mbox{ or } 1 \leq i=j \leq m\})+ n-m.
\end{eqnarray*}
Also, since $e_{m+1},\dots,e_n$ appear only in $\pusai(m+1,m+1),\dots,\pusai(n,n)$,
we have 
\begin{eqnarray*}
s_{g,k}(G)&=&\dim_K (\mathrm{span}\{\pusai(i,j):\{i,j\} \in E(G) \mbox{ or } 1 \leq i=j \leq n\})\\
&=&\dim_K (\mathrm{span}\{\pusai(i,j):\{i,j\} \in E(G) \mbox{ or } 1 \leq i=j \leq m\}) + n-m.
\end{eqnarray*}
For a generic matrix $g \in \GL$ and for $ 2 \leq k \leq m$,
since $f_k(i,j)$ is coincident with $\pusai (i,j)$,
(\ref{sono4}) and (\ref{ppc2}) say that $m_k(\GIN(G))=s_{g,k}(G)$.
Also, for $k >m$,
Lemma \ref{tyoutenn} says that $m_{\leq k} (\GIN(G))=m_{\leq m} (\GIN(G))=|E(G)|+n$.
Also, since $s_{g,k+1}(G) \geq s_{g,k}(G)$ and $s_{g,k}(G) \leq |E(G)|+n$ hold for all $k$,
we have 
$$|E(G)|+n \geq s_{g,k}(G) \geq s_{g,m}(G)= m_{\leq m} (\GIN(G))=|E(G)|+n \quad\mbox{ for } k >m.$$
Thus we have $s_{g,k}(G)=m_{\leq k} (\GIN(G))=|E(G)|+n$ for $k >m$.
In particular, by Lemma \ref{symmetric},
 we have $s_{g,k+1}(G)-n=m_{\leq k } \dels G$.
\end{proof}

\section{Algebraic shifting of bipartite graphs}
In this section we study algebraic shifting of bipartite graphs.

\begin{lemma} \label{binomial}
For any positive integer $n>0$,
there exist unique integers $h(n)$ and $\alpha(n)$ such that
$$n= {h(n) \choose 2} +\alpha(n) \ \ \ \ \ \mbox{ and }\ \  h(n) \geq \alpha(n) > 0.$$
\end{lemma}
 
\begin{proof} 
Let $h(n)$ be the maximal integer which satisfies ${ h(n) \choose 2} <n$.
Set $\alpha(n)=n- {h(n) \choose 2}>0$.
Since $n \leq {h(n)+1 \choose 2}={h(n) \choose 2}+h(n)$,
we have $\alpha(n) \leq h(n)$.
Then $h(n)$ and $\alpha(n)$ satisfy $n={h(n) \choose 2} +\alpha(n)$ and $h(n) \geq \alpha(n) > 0$.

Conversely, if $h(n)$ and $\alpha(n)$ satisfy $n={h(n) \choose 2} +\alpha(n)$ and $h(n) \geq \alpha(n) > 0$,
then $h(n)$ satisfies ${h(n)+1 \choose 2} \geq n > {h(n) \choose 2}$
and such an integer $h(n)$ is uniquely determined.
Thus  $\alpha(n)=n- {h(n) \choose 2}$ is also uniquely determined.
\end{proof} 

We call the  above representation
% of an integer $n >0$
the \textit{binomial form} of $n$.

For bipartite graphs, we have the following relation between exterior algebraic shifting and symmetric algebraic shifting.

\begin{lemma}\label{relation}
Let $G$ be a bipartite graph on $[n]$.
Then, for $1 \leq k \leq n-1$, one has 
\begin{eqnarray}
\ \ \ \ \ \ m_{\leq k+1}(\dele G) \geq m_{\leq k}(\dels{G})\geq m_{\leq k+1} (\dele{G})-n+\min\{ {{k+2}\choose 2} , n\} \label{bound}
\end{eqnarray}
\end{lemma}

\begin{proof}
Assume that $G$ is a bipartite graph on the vertex set $A\cup B =[n]$
with $E(G)\cap \{ E(K_A) \cup E(K_B)\} = \emptyset$,
where $K_A$ is the complete graph on the vertex set $A$ and $K_B$ is the complete graph on the vertex set $B$.
Let $V$ be the $n$-dimensional $K$-vector space with basis $e_1,\dots ,e_n$,
$V_A$ the subspace of $V$ spanned by $e_i$ with $i \in A$
and $V_B$ the subspace of $V$ spanned by $e_i$ with $i \in B$. 

Fix an integer $1 \leq k \leq n-1$.
Let 
$g=(x_{ij}) \in \GL$.
First, we will show
\begin{eqnarray}
\dim_K \{\pusai(G)\}=r_{g,k}(G). \label{srm3}
\end{eqnarray}
For integers $a\in A$ and $b \in B$,
let
\begin{eqnarray*}
\tuua(a,b) =(x_{1b}e_a,x_{2b}e_a,\dots,x_{kb}e_a) \ \ 
% \in \bigoplus_{i=1}^k V_A \\
\mbox{and} \ \ \tuub(a,b)=(x_{1a}e_b,x_{2a}e_b,\dots,x_{ka}e_b)
% \in \bigoplus_{i=1}^k V_B.
\end{eqnarray*}
Then
$$\phai (a,b) = \tuua (a,b) - \tuub (a,b) \ \ \ 
\mbox{ and }
\ \ \ \pusai (a,b) = \tuua (a,b) + \tuub (a,b).$$
Let $\Phi$ be the automorphism of the $K$-vector space $V$ defined by
$\Phi(e_i) =e_i$ if $i \in A$, and $\Phi(e_i) = -e_i$ if $i \in B$.
Let $\Phi_k:\bigoplus_{i=1}^k V \to \bigoplus_{i=1}^k V$ be
the automorphism of  $\bigoplus_{i=1}^k V$ defined by
$\Phi_k(v_1,\dots,v_k)=(\Phi(v_1),\dots,\Phi(v_k))$.
Since
$\Phi_k ( \pusai (a,b))=\phai(a,b)$,
we have $ \Phi_k ( \pusai(G)) =\phai (G)$.
Since $\Phi_k$ is an automorphism, we have 
\begin{eqnarray}
\dim_K(\pusai(G))=\dim_K (\phai(G))=r_{g,k}(G).\label{go}
\end{eqnarray}
Then the definition of $s_{g,k}(G)$ says that
$$s_{g,k}(G)= \dim_K (\pusai(G^+)) \leq \dim_K (\pusai (G)) +n =r_{g,k}(G)+n.$$
Hence, by Lemmas \ref{exdet} and \ref{GIN}, for a generic matrix $g \in \GL$, we have
\begin{eqnarray*}
m_{\leq k} (\dels G) =s_{g,k+1}(G) -n \leq r_{g,k+1}(G)=m_{\leq k+1}(\dele G).
\end{eqnarray*}
This is the first inequality of (\ref{bound}).\bigskip

Next, we consider the second equality of (\ref{bound}).
We use the following idea appeared in \cite[Theorem 6.1]{KA-h}.
For each $1 \leq i \leq n$,
define a linear map $f_i^A:V \to K$ by
\begin{eqnarray*}
f_i^A(e_k)=
\left\{
\begin{array}{l}
x_{ik}, \ \ \ \mbox{ if } k \in A,\\
0, \ \ \ \ \ \ \mbox{otherwise.}
\end{array}
\right.
\end{eqnarray*}
Also we define $f_i^B: V \to K$ by the same way.
For integers $1 \leq i \leq j \leq k$,
let $F_k^{(i,j)}: \bigoplus_{l=1}^k V \to K$ be the linear map defined by
$$F_k^{(i,j)}(v)= f_i^A(v_j) -f_j^B(v_i) \ \ \ \ \ \ \ \mbox{ for } \ v=(v_1,v_2,\dots,v_k) \in \bigoplus_{i=1}^k V.$$
Consider the linear map $F_{g,k}:\ \pusai (G^+) \to \bigoplus_{1\leq i\leq j \leq k} K$
defined by
$$
F_{g,k}(v)=(F_k^{(1,1)}(v), F_k^{(1,2)}(v),F_k^{(2,2)}(v),F_k^{(1,3)}(v),F_k^{(2,3)}(v),F_k^{(3,3)}(v),\dots,F_k^{(k,k)}(v)).
$$
Then, for any $\{a,b\} \in E(G)$ and for all integers $1 \leq i \leq j \le k$,
we have
$$
F_k^{(i,j)} ( \pusai ( a,b)) =x_{jb} x_{ia} -x_{ia} x_{jb} =0 .
$$
Thus we have $\pusai(G) \subset \ker(F_{g,k} )$.
Hence, by (\ref{go}), we have
$\dim_K (\ker(F_{g,k})) \geq r_{g,k}(G)$.
Then, by Lemmas \ref{exdet} and \ref{GIN},
for a generic matrix $g \in \GL$,
we have
\begin{eqnarray}
m_{\leq k-1}(\dels G)&=&\dim_K\{\pusai(G^+)\}-n \nonumber \\
&=& \dim_K \{ \ker (F_{g,k})\} + \dim_K \{ \mathrm{Im}(F_{g,k})\} -n  \label{srm5} \\
&\geq& m_{\leq k} (\dele G)+ \dim_K \{ \mathrm{Im}(F_{g,k})\}-n.\nonumber
\end{eqnarray}
We claim
$\dim_K ( \mathrm{Im}(F_{g,k})) \geq \min \{ { k+1 \choose 2},n\}$ for a generic matrix $g \in \GL$.

Since $\pusai(G) \subset \ker (F_{g,k})$, it follows that $\mathrm{Im}(F_{g,k})=F_{g,k}(\pusai([n]))$.
Also, for each $1 \leq t \leq n$, the vector $F_{g,k}(\pusai (t,t))$ is of the form
\begin{eqnarray*}
F_{g,k}(\pusai(t,t))=
\left \{
\begin{array}{l}
(x_{1t}x_{1t},x_{1t}x_{2t},x_{2t}x_{2t},x_{1t}x_{3t},x_{2t}x_{3t},x_{3t}x_{3t},\dots,x_{kt}x_{kt}), \ \ \hspace{1.5pt}\mbox{ if } t \in A,\\
-(x_{1t}x_{1t},x_{1t}x_{2t},x_{2t}x_{2t},x_{1t}x_{3t},x_{2t}x_{3t},x_{3t}x_{3t},\dots,x_{kt}x_{kt}), \mbox{ if } t \in B.
\end{array}
\right.
\end{eqnarray*}
Let $M_{g,k}= (x_{it}x_{jt})_{1 \leq t \leq n,\ 1 \leq i \leq j \leq k}$ be the $n \times {k+1 \choose 2}$ matrix
whose $(t,(i,j))$-th entry is $x_{it}x_{jt}$.
Since each $t$-th low vector of $M_{g,k}$ is coincident with $F_{g,k}( \pusai (t,t))$,
we have 
$$ \dim_K (\mathrm{Im}(F_{g,k}))=\dim_K\{ F_{g,k}(\pusai ([n]))\}= \mathrm{rank}(M_{g,k}).$$

Let $l = \min\{ { k+1 \choose 2},n\}$.
For $1 \leq t \leq n$,
let $t={h(t) \choose 2} +\alpha(t)$ be the binomial form of $t$ defined in Lemma \ref{binomial}
and let $\pi(t)=(\alpha(t),h(t))$.
Note that $\pi(1)=(1,1),\pi(2)=(1,2),\pi(3)=(2,2),\pi(4)=(1,3),\dots$.
Let $\widetilde M_{g,k}$ be the $l \times l$ submatrix of $M_{g,k}$ whose lows are indexed by $1,2,\dots,l$ and whose
columns are indexed by $\pi(1),\pi(2),\dots,\pi(l)$.
Then the product of $(t,\pi(t))$-th entries of $\widetilde M_{g,k}$ for $t=1,2,\dots,l$ is 
$$u=x_{\alpha(1)1}x_{h(1)1}x_{\alpha(2)2}x_{h(2)2} \cdots x_{\alpha(l)l}x_{\alpha(l)l}.$$
We regard the determinant $\det (\widetilde M_{g,k})$ of $\widetilde M_{g,k}$ as a polynomial in $K[x_{ij}]_{1\leq i,j\leq n}$.
Then the definition of $\widetilde M_{g,k}$ says that the coefficients of the monomial 
$u$ in $\det(\widetilde M_{g,k}) \in K[x_{ij}]_{1\leq i,j\leq n}$ is $1$.
Thus $\det(\widetilde M_{g,k})$ is not zero as a polynomial,
 and therefore we \vspace{1.5pt} have $\mathrm{rank}(\widetilde M_{g,k})=l$ for a generic matrix $g \in \GL$.
Hence we have
$$\dim_K(\mathrm{Im}(F_{g,k}))= \mathrm{rank}(M_{g,k}) \geq \mathrm{rank}(\widetilde M_{g,k})= \min \{ { k+1 \choose 2},n\},$$
for a generic matrix $g \in \GL$ as desired.
Thus the assertion follows from (\ref{srm5}).
\end{proof}

Exterior algebraic shifting of complete bipartite graphs are determined as follows

\begin{lemma}[{\cite[Theorem 6.1]{KA-h}}]\label{bipar}
Let $a \geq b >0$ be positive integers and $n=a+b$.
Let $K_{a,b}$ be the complete bipartite graph of size $a,b$.
Then one has
\begin{eqnarray*}
m_{\leq k}(\dele{K_{a,b}})=
\left \{
\begin{array}{l}
kn-k^2,\ \mathrm{if} \ k\leq b, \\
ab,\ \ \  \ \ \ \ \hspace{3pt}        \mathrm{if} \ k>b .
\end{array}
\right.
\end{eqnarray*}
\end{lemma}

Before considering symmetric algebraic shifting of complete bipartite graphs,
we recall the following property.
Let $G$ be a graph on $[n]$.
For a vertex $v \in [n]$ of $G$,
write $\deg_G(v) =|\{t\in [n]: \{v,t\} \in E(G)\}|$
and write $G-\{v\}$ for the induced subgraph of $G$ on $[n] \setminus \{v\}$.

\begin{lemma} \label{hyperconnect}
Let $G$ be a graph on $[n]$ and $v\in [n]$.
If  $ \{k+1,k+2\}\not\in \dels {G-\{v\}}$ and $\deg_G(v)\leq k$
then one has $\{k+1,k+2\}\not\in \dels G$. 
\end{lemma}

In \cite[Lemma 4.3]{KA-h}, Kalai proved the above lemma for exterior algebraic shifting.
However, the proof for symmetric algebraic shifting is the same.

\begin{theorem}\label{c-bipar}
Let $a \geq b >0$ be positive integers and $n=a+b$.
Let $n= {h(n) \choose 2} +\alpha(n)$ be the binomial from of $n$. 
Then
\begin{eqnarray*}
m_{\leq k}(\dels{K_{a,b}})=
\left \{
\begin{array}{l} \smallskip
{n \choose 2}-{n-k \choose 2},\ \ \ \ \ \ \  \ \ \ \  \hspace{2pt} \mathrm{if} \ k\leq b-1 \ \mathrm{and} \  k\leq h(n)-2, \\ \smallskip
(k+1)n-(k+1)^2,\ \  \mathrm{if} \ k\leq b-1 \ \mathrm{and} \  k> h(n)-2, \\
ab,\ \ \ \hspace{83pt}       \mathrm{if} \ k>b-1.
\end{array}
\right.
\end{eqnarray*}
\end{theorem}
\begin{proof}
If $k >h(n)-2$, then $\min \{ {{k+2} \choose 2},n\}=n$.
Then Lemma \ref{relation} says
$$m_{\leq k} (\dels {K_{a,b}}) =m_{\leq k+1} (\dele {K_{a,b}}).$$
Thus, by Lemmas \ref{relation} and \ref{bipar}, Theorem \ref{c-bipar} is true for $k > h(n) -2$.
\medskip

Next, we assume $k \leq h(n)-2$ and $k \leq b-1$.
Then Lemmas \ref{relation} and \ref{bipar} says that
\begin{eqnarray}
m_{\leq k}(\dels {K_{a,b}}) &\geq& m_{\leq k+1}(\dele {K_{a,b}}) + {{k+2} \choose 2}-n  \label{srmt} \\
&=& kn+(k+1)^2 + {{k+2} \choose 2}. \nonumber 
\end{eqnarray}
Also, a routine computation says that
$kn+(k+1)^2 + {{k+2} \choose 2} = {n \choose 2} - {n-k \choose 2}$.
On the other hand,
for an arbitrary graph $G$ on $[n]$,
one has
\begin{eqnarray}
m_{\leq k}(G) \leq (n-1)+(n-2) +\cdots + (n-k)={n \choose 2} - {n-k \choose 2}.\label{srms}
\end{eqnarray}
Then inequalities (\ref{srmt}) and (\ref{srms}) say that
$m_{\leq k} (\dels {K_{a,b}}) ={n \choose 2} - {n-k \choose 2}$ as desired.
\medskip

At last, we consider the case $k \leq h(n)-2$ and $k > b-1$.
We use induction on $a+b$.
The assertion is obvious if $a+b=2$.
Assume $a+b >2$ and consider $K_{a-1,b}$.
Then the assumption of induction says that
$$m_{\leq k}(\dels {K_{a-1,b}}) =(a-1)b=|E(\dels{K_{a-1,b}})| \ \ \ \ \ \ \mbox{ for } \ k > \min\{a-1,b\}-1.$$
Thus, in particular, we have $\{b+1,b+2\} \not\in E(\dels {K_{a-1,b}})$.
Then Lemma \ref{hyperconnect} says that $\{b+1,b+2\} \not \in E( \dels {K_{a,b}})$.
Since $\dels {K_{a,b}}$ is shifted, any $\{i,j\} \subset [n]$ with $\min\{i,j\} \geq b+1$ does not belong to $E(\dels {K_{a,b}})$.
Thus we have 
$$m_{\leq k}(\dels {K_{a,b}})=|E(\dels {K_{a,b}})|=ab \ \ \ \ \ \ \mbox{ for } \ k > b-1,$$
 as desired.
\end{proof}

\begin{cor}\label{coro}
Let $G$ be a bipartite graph on $[n]$.
Let $n= {h(n) \choose 2} +\alpha(n)$ be the binomial from of $n$. 
If $\{h(n),h(n)+1\} \in \dele{G}$, then $\dele{G}\ne \dels{G}$.
\end{cor}
\begin{proof}
Since $n< { h(n) +1 \choose 2}$, Lemma \ref{relation} says that
$m_{\leq h(n)-1}(\dels G) = m_{\leq h(n)}(\dele{G})$.
On the other hand,
since $\{h(n),h(n)+1\} \in \dele{G}$,
we have $m_{\leq  h(n)}(\dele{G})>m_{\leq  h(n)-1}(\dele{G})$.
Thus we have  $m_{\leq  h(n)-1}(\dels{G})>m_{\leq  h(n)-1}(\dele{G})$.
\end{proof}

\begin{exam}
By Lemma \ref{bipar} and Theorem \ref{c-bipar},
we can easily compare \dele{K_{a,b}} and \dels{K_{a,b}}.
For example,
\begin{eqnarray*}
\dele{K_{6,6}}=
\left\{
\begin{array}{l}
\{1,2\},\dots,\{1,12\},\{2,3\},\dots,\{2,11\},\{3,4\},\dots,\{3,10\},\\
\{4,5\},\dots,\{4,9\},\{5,6\},\dots,\{5,8\},\{6,7\}\}
\end{array}
\right\}.
\end{eqnarray*}
and 
\begin{eqnarray*}
\dels{K_{6,6}}=
\left\{
\begin{array}{l}
\{1,2\},\dots,\{1,12\},\{2,3\},\dots,\{2,12\},\{3,4\},\dots,\{3,12\},\\
\{4,5\},,\dots\{4,9\},\{5,6\}\}.
\end{array}
\right\}.
\end{eqnarray*}

Also, Corollary \ref{coro} implies that,
for many of the bipartite graphs, their exterior algebraic shifted graph and symmetric algebraic shifted graph
are different.
If  $G$ is a shifted graph with $|E(G)|>{n \choose 2} - {n-h(n)+1 \choose 2}$,
then $\{h(n),h(n)+1\}$ must be contained in $E(G)$.
Thus if $G$ is a bipartite graph on $[n]$ with
 $|E(G)|>{n \choose 2} - {n-h(n)+1 \choose 2}$, then one has $\dele G \ne \dels G$.

Assume that $G$ is a subgraph of $K_{n,n}$.
Then
\begin{eqnarray*}
{2n \choose 2} - {2n-h(2n)+1 \choose 2}=  2n(h(2n)-1) -{h(2n) \choose 2}. \label{srmhosi1}
\end{eqnarray*}
On the other hand, the definition of binomial form says that ${h(2n) \choose 2} <2n$, and therefore $(h(2n)-1)^2 < 4n$.
Thus $|E(K_{n,n})|=n^2 \gg 2n(h(2n)-1) -{h(2n) \choose 2}$ for $n \gg 0$.
This fact says that
exterior algebraic shifting and symmetric algebraic shifting of most of the subgraphs of $K_{n,n}$ are
different if $n \gg 0$.
\end{exam}

\section{Chordal graphs and graded Betti numbers}

Let $K$ be a field of characteristic $0$ and
$R=K[x_1,x_2,\dots,x_n]$  the polynomial ring.
The graded Betti numbers $\beta_{ij}(I)$ of a homogeneous ideal $I\subset R$
are the integers $\beta_{ij}(I)=\dim_K (\mathrm{Tor}_i(I,K)_j)$.
In other words,
$\beta_{ij}(I)$ appears in the minimal graded free resolution
\[
0 \longrightarrow \bigoplus_{j} R(-j)^{\beta_{hj}(I)}
\longrightarrow
\cdots
\longrightarrow
\bigoplus_{j} R (-j)^{\beta_{1j}(I)}
\longrightarrow
\bigoplus_{j} R (-j)^{\beta_{0j}(I)}
\longrightarrow I
\longrightarrow 0
\]
of $I$ over $R$.

Let $G$ be a graph on $[n]$.
Write $I(G)$ for the ideal generated by all squarefree monomials $x_i x_j$ with $\ij \not\in E(G)$.
If $G$ is shifted, the second graded Betti numbers of $I(G)$ are given by the formula (\cite[Corollary 3.6]{Her})
\begin{eqnarray*}
\beta_{i i+2}( I(G)) &=&
\sum_{s=1}^n |\{x_s x_t\in I_G: s<t\}| { n+1-s-1 \choose i}\\
&=& \sum_{s=1}^n \{n-s - m_s(G)\} { n+1-s-1 \choose i}.
\end{eqnarray*}

The above formula implies that
if $G$ and $G'$ are shifted graphs on $[n]$ which satisfy 
$\beta_{ii+2}(I(G))=\beta_{ii+2}(I(G'))$ for all $i \geq 0$,
then one has $m_k(G) = m_k(G')$ for all $k=1,2,\dots,n$.
Thus the relation (\ref{tondemo}) before Lemma \ref{exdet} guarantees

\begin{lemma}\label{kimeru}
Let $G$ and $G'$ be shifted graphs on $[n]$.
The followings are equivalent.
\begin{itemize}
\item[(i)] $\beta_{ii+2}(I(G))=\beta_{ii+2}(I(G)) $ for all $i\geq 0$;
\item[(ii)] $m_{\leq k}(G)=m_{\leq k}(G')$ for all $k$;
\item[(iii)] $G=G'$.
\end{itemize}
\end{lemma}

Let $I$ be a homogeneous ideal generated in degree $d$.
We say that $I$ has a \textit{linear resolution} if $\beta_{ii+j}(I)=0$ for all $j\ne d$.

Let $G$ be a graph and $(\{r_1,r_2\},\dots,\{r_{s},r_{s+1}\})$ a cycle in $G$, where $r_{s+1}=r_1$.
An edge $\ij \in E(G)$ is called a \textit{chord} of the cycle $(\{r_1,r_2\},\dots,\{r_{s},r_{s+1}\})$
if $i,j\in \{r_1,r_2,\dots,r_s\}$ and $\ij \ne \{r_i,r_{i+1}\}$ for all $i$.
We call $G$ \textit{chordal} if every cycle in $G$ of length $>3$ has a chord.
Note that shifted graphs are chordal.
The following fact is known.

\begin{lemma}\label{codal}
Let $G$ be a graph on $[n]$.
The followings are equivalent.
\begin{itemize}
\item[(i)] $G$ is chordal;
\item[(ii)] $I(G)$ has a linear resolution;
\item[(iii)] $\beta_{ii+2}(I(G))=\beta_{ii+2}(I(\dele G))=\beta_{ii+2}(I(\dels{G}))$ for all $i\geq 0$.
\end{itemize}
\end{lemma}

\begin{proof}
(i) $\Leftrightarrow$ (ii) was proved in \cite{chordal}.
We will show (ii) $\Leftrightarrow$ (iii). 
$I(G)$  is the ideal generated by all monomials in $I_G$ of degree $2$.
Thus, it follows from \cite[Theorem 2.1 and Corollary 2.3]{stablenumber} that
$I(G)$ has a linear resolution
if and only if
$\beta_{ii+2}(\mathrm{gin}(I_G))=\beta_{ii+2}(I_G)$ and $\beta_{ii+2} (I_{\dele G}) =\beta_{ii+2}(I_G)$
for all $i$.
On the other hand, it follows from \cite[Lemma 8.18]{Her} that
$\beta_{ij}(\mathrm{gin}(I_G))=\beta_{ij}( I_{\dels G})$ for all $i,j$.
%Also, it follows from \cite[Lemma 1.2]{componentwise} that
Also, since
% it is easy to see that
$\beta_{i i+2}(I_G)= \beta_{ii+2} (I(G)) $, $\beta_{i i+2}(I_{\dels G})= \beta_{ii+2} (I(\dels G))$
and $\beta_{i i+2}(I_{\dele G})= \beta_{ii+2} (I(\dele G))$
for all $i \geq 0$ (see \cite[Lemma 1.2]{componentwise}),
it follows that conditions (ii) and (iii) are equivalent. 
\end{proof}

Lemma \ref{kimeru} together with Lemma \ref{codal}
immediately implies

\begin{lemma}\label{onazi}
If $G$ is a chordal graph,
then one has $\dele G =\dels G$.
\end{lemma}

\begin{exam}
In general, algebraic shifting is defined for simplicial complexes
(see \cite{Her} or \cite{KA}).
Lemmas \ref{kimeru} and \ref{codal} say that
if the Stanley--Reisner ideal $I_\Gamma$ of a simplicial complex $\Gamma$ are generated in degree $2$ and if $I_\Gamma$ has a linear resolution,
then one has $\dele \Gamma = \dels \Gamma $.
However,
there exists a simplicial complex $\Gamma$ such that
$I_\Gamma$ has a linear resolution, but  $\dele \Gamma \ne \dels \Gamma$,
if $I_\Gamma$ is generated in degree $d \geq 4$.

Let $\Gamma$ be the simplicial complex generated by
\begin{eqnarray*}
\Big\{
\{3,4,5,6\},\{2,4,5,6\},\{1,4,5,6\},\{1,2,3,4\},\{1,2,3,5\},\{1,2,3,6\}
\Big\}
\end{eqnarray*}
Note that $\Gamma$ is the Alexander dual of $K_{3,3}$.
Then, $I_\Gamma$ has a liner resolution and generated in degree $4$.
The $4$-skeleton of $\dele \Gamma$ is
$$
\Big\{
\{1,2,3,4\},\{1,2,3,5\},\{1,2,3,6\},\{1,2,4,5\},\{1,2,4,6\},\{1,3,4,5\}
\Big\},
$$
however, the $4$-skeleton of $\dels \Gamma$ is 
$$ 
\Big\{
\{1,2,3,4\},\{1,2,3,5\},\{1,2,3,6\},\{1,2,4,5\},\{1,2,4,6\},\{1,2,5,6\}
\Big\}.$$
( The form of $\dele \Gamma$ follows from
the relation between exterior algebraic shifting and Alexander dual \cite[pp. 137]{KA}.
For $\dels \Gamma$, we compute $\gin (I_\Gamma)$ by using \cocoa)

\end{exam}

\section{Combinatorial shifting of chordal graphs}

In this section, we study combinatorial shifting of chordal graphs.
First, we will determine all chordal graphs whose combinatorial shifted graph is unique.
Second, we will show that, for any chordal graph $G$,
there exist a combinatorial shifted graph $\delc G$ such that $\delc G = \dele G$.
In particular, we will give an algorithm to compute the combinatorial shifted graph $\delc G$ of a chordal graph $G$
with $\delc G =\dele G$.

Combinatorial shifting was introduced by Erd\"os, Ko and Rado \cite{EKR}.
Let $\Gamma$ be a collection of subsets of $[n]$.
For integers $1 \leq i < j \leq n$,
write $\mathrm{Shift}_{ij}(G)$
for the collection of subsets of $[n]$
whose elements are $C_{ij}^\Gamma(S) \subset [n]$,
where $S \in G$ and where
\begin{eqnarray*}
C_{ij}^\Gamma(S) =
\left\{
\begin{array}{l}
(S \setminus \{ j \}) \cup \{ i \},
\, \, \, \, \, \mbox{if} \, \, \,
j \in S, \, \, \,
i \not\in S \, \, \,
\mbox{and} \, \, \,
(S \setminus \{ j \}) \cup \{ i \} \not\in \Gamma,
\\
\hspace{0cm} S, \hspace{2.6cm}
\mbox{otherwise.}
\end{array}
\right.
\end{eqnarray*}
Let $G$ be a graph on $[n]$.
Write $\shif$ for the finite graph on $[n]$ with the edge set $E( \shif) =\shift_{ij}(E(G))$.
It follows from, e.g., \cite[Corollary 8.6]{Her}
that there exists a finite sequence of pairs of integers
$(i_1, j_1), (i_2, j_2), \ldots, (i_q, j_q)$
with each $1 \leq i_k < j_k \leq n$ such that
\[
\mathrm{Shift}_{i_qj_q} (\mathrm{Shift}_{i_{q-1}j_{q-1}}( \cdots
(\mathrm{Shift}_{i_1j_1}(G) )\cdots ))
\]
is shifted.  Such a shifted graph is called
a {\em combinatorial shifted graph} of $G$,
and will be denoted by $\Delta^c(G)$.
Combinatorial shifted graphs are easily computed,
however,
they are not necessarily unique.

Let $\sigma :[n] \to [n]$ be a permutation of $[n]$.
Write $\sigma(G)$ for the graph on $[n]$
with the edge set $E(\sigma(G))= \{ \{\sigma(s),\sigma(t)\}: \{s,t\} \in E(G)\}$.
We say that a graph $G'$ is isomorphic to $G$ if there exists a permutation $\sigma:[n] \to [n]$
such that $\sigma(G')=G$.
Write $\sigma_{ij}: \nn \to \nn$ for the transposition of $i,j\in \nn$.

\begin{lemma} \label{ape}
Let $G$ and  $H$ be graphs on $[n]$.
If $G$ is isomorphic to $H$,
then any combinatorial shifted graph of $G$ is equal to some
combinatorial shifted graph of $H$.
%such that $\delc G =\delc H$.
\end{lemma}

\begin{proof}
We may assume $H=\sigma_{ij}(G)$ for some $i,j\in \nn$.
Let 
$$r(G)=|\{ \ij \subset [n]:
\ij \not \in E(G) \mbox{ and } \{ i',j'\} \in E(G) \mbox{ for some } i' \geq i, j' \geq j \} |.$$
We use induction on $r(G)$.
If $r(G)=0$, then $G$ is shifted and the only combinatorial shifted graph of $G$ is $G$ itself.
On the other hand,
for any integers $1 \leq p < q \leq n$,
the definition of $\shift_{pq}(G)$ says that
\begin{eqnarray}
\hspace{30pt}
E(\shift_{pq}(G))=
\left\{
\begin{array}{l}
\ \{ \{s,t\} \in E(G): \{s,t\}=\{p,q\} \mbox{ or } p,q \not\in \{s,t\}\}\smallskip\\
\cup\ \{ \{p,t\}: t \ne q,\ \{p,t\} \in E(G) \mbox{ or } \{ q,t\} \in E(G)\}\}\smallskip\\
\cup\ \{ \{q,t\}: t \ne p,\ \{p,t\} \in E(G) \mbox{ and } \{ q,t\} \in E(G)\}\}.
\end{array}
\right. \label{ml4-1} 
\end{eqnarray}
The above form together with (\ref{tondemo}) implies $G=\shift_{ij}(\sigma_{ij}(G))=\shift_{ij}(H)=\delc H$.

Next, we assume $r(G)>0$.
Let $\delc G= \shift_{i_lj_l}(\cdots(\shift_{i_1j_1}(G))\cdots)$.
Assume $G \ne \shift_{i_1j_1}(G)$.
Then we have $r( \shift_{i_1j_1}(G))<r(G)$.

If $\sigma_{ij}(i_1)<\sigma_{ij}(j_1)$, then, by (\ref{ml4-1}),
we have $\shift_{ \sigma_{ij}(i_1)\sigma_{ij}(j_1)}(H)= \sigma_{ij}\{ \shift_{i_1j_1}(G)\}$.
Let $H_1=\shift_{ \sigma_{ij}(i_1)\sigma_{ij}(j_1)}(H)$.
Then $\shift_{i_1j_1}(G)$ is isomorphic to $H_1$.
Thus the induction hypothesis says that
there exists a combinatorial shifted graph $\delc {H_1}$ of $H_1$
such that $\delc G = \delc {H_1}$. Since $\delc {H_1}$ is a combinatorial shifted graph of $H$,
the claim follows.

If $\sigma_{ij}(i_1)>\sigma_{ij}(j_1)$,
then, by (\ref{ml4-1}),
we have
$\sigma_{\sigma_{ij}(j_1) \sigma_{ij}(i_1)}(\shift_{\sigma_{ij}(j_1) \sigma_{ij}(i_1)}(H))= \sigma_{ij} (\shift_{i_1j_1}(G)).$
Then $\shift_{\sigma_{ij}(j_1) \sigma_{ij}(i_1)}(H)$ is isomorphic to $\shift_{i_1 j_1}(G)$.
Thus, by the same way as the case $\sigma_{ij}(i_1)>\sigma_{ij}(j_1)$,
% the assertion follows.
there is $\delc H$ such that $\delc G = \delc H$.
\end{proof}

Let $G$ be a graph on $\nn$.
We call $v\in \nn$  a \textit{star vertex of $G$}
if $\{u,v\}\in E(G)$ for all $u\in [n]\setminus \{v\}$ with $\deg_G(u) > 0$. 
Also, we call $v \in [n]$ an \textit{isolated vertex of} $G$ if $\deg_G(v)=0$. 
Let $\overline{G-\{v\}}$ be the induced subgraph of $\sigma_{1v}(G)$ on $[2,n]=\{2,\dots,n\}$.

\begin{lemma}\label{central}
Let $G$ be a  graph on $[n]$ and
$v \in [n]$ a star vertex of $G$.
Then,
for every combinatorial shifted graph $\delc G$ of $G$,
there exists a combinatorial shifted graph $\delc {\overline {G-\{v\}}}$ of $\overline {G-\{v\}}$ on $[2,n]$
such that
\begin{eqnarray}
E(\delc G) = E(\delc {\overline {G-\{v\}}}) \cup \{ \{1,2\},\dots,\{1,1+\deg_G(v)\} \}.\label{ml44}
\end{eqnarray}
Conversely, for any combinatorial shifted graph $\delc {\overline {G-\{v\}}}$ of $\overline {G-\{v\}}$
there exists a combinatorial shifted graph $\delc G$ which satisfies (\ref{ml44}).
\end{lemma}

\begin{proof}
By Lemma \ref{ape}, we may assume $v=1$.
Let $H$ be the subgraph of $G$ on $[n]$ with $E(H)=\{\{s,t\}\in E(G): 1 \in \{s,t\} \}$.
For integers $1 < i < j \leq n$, we have
\begin{eqnarray}
E(\shift_{ij}(G))=E(\shift_{ij}(\overline {G-\{1\}})) \cup E(\shift_{ij}(H)). \label{also}
\end{eqnarray}
Note that the definition of $\shif$ says that $1\in \nn$ is a star vertex of $\shift_{ij}(G)$.
Also,
since $1$ is a star vertex of $G$, we have $\shift_{1j}(G)=G$ for any $1<j \leq n$.
These facts say that,
for any combinatorial shifted graph $\delc G$,
there exist combinatorial shifted graphs $\delc {\overline {G-\{1\}}}$ and $\delc H$ such that
\begin{eqnarray}
E(\delc G)= E(\delc {\overline { G-\{1\} } }) \cup E(\delc H). \label{zyuusi}
\end{eqnarray}
Conversely, (\ref{also}) also says that, for any combinatorial shifted graph $\delc { \overline{G-\{1\}}}$,
there is a combinatorial shifted graph $\delc G$ and $\delc H$ which satisfy (\ref{zyuusi}).
Since every edge $\{s,t\} \in E(\shift_{ij}(H))$ contains $1$ for all $i,j$,
the graph $H$ has the unique combinatorial shifted graph $\delc H$ with
$E(\delc H)=\{ \{1,2\},\dots,\{1,1+\deg_G(1)\}\}$.
Thus the assertion follows.
\end{proof}

We will determine all chordal graphs whose combinatorial shifted graph is unique.
First, we give an easy example.

\begin{lemma}\label{disj}
Let $A \subset [n]$ and $B \subset [n]$ are subsets with $A \cap B = \emptyset$.
Let $G=K_A\cup K_B$ be the graph with $E(G)= E(K_A) \cup E(K_B)$.
Then $G$ has the unique combinatorial shifted graph.
\end{lemma}

\begin{proof}
Let $\delc G =\shift_{i_pj_p} (\cdots (\shift_{i_1j_1}(G))\cdots )$.
If $\{ i_1,j_1\} \subset A$ or $\{ i_1,j_1\} \subset B$,
then $\shift_{i_1j_1} (G)=G$.
Thus we may assume $i_1 \in A$ and $j_1 \in B$.
Let $H= \shift_{i_1j_1}(G)$.
Then, by the definition of $\shift$, the vertex $i_1$ is a star vertex of $\shift_{i_1j_1}(G)$.
Also, we have  $E(H-\{i_1\})=E(K_{A \setminus \{i_1\}})\cup E(K_{B \setminus \{j_1\}})$.
Then Lemma \ref{central}
says  that there exist a combinatorial shifted graph $\overline{\delc {H-\{i_1\}}}$ such that
$$E(\delc G) = E(\delc {\overline {H-\{i_1\}}}) \cup \{ \{1,2\}, \dots,\{1,|A|+|B|\} \}.$$
Arguing inductively,
the claim follows.
% it follows that $G$ has the unique combinatorial shifted graph.
\end{proof}

\begin{rem}
The complementary graph of a complete bipartite graph 
satisfy the condition of Lemma \ref{disj}.
Since the complementary graph of a combinatorial shifted graph $\delc G$ of $G$
is equal to some combinatorial shifted graph of the complementally graph of $G$,
complete bipartite graphs have the unique combinatorial shifted graph.
Also, it is easy to see that this unique combinatorial shifted graph $\delc{K_{a,b}}$
is equal to $\dele {K_{a,b}}$.
\end{rem}

Let $G$ be a graph on $[n]$.
We call $\shift_{ij}$ an \textit{edge shift of} $G$ if $\ij \in E(G)$.
Also, we call $\shift_{ij}$ a \textit{disjoint shift of} $G$ if
there is no path from $i$ to $j$ in $G$.

\begin{lemma}\label{connect}
Let $G$ be a connected graph on $[n]$.
Then there exist a sequence 
$(i_1, j_1), (i_2, j_2), \ldots, (i_q, j_q)$
of pairs of integers with each $1 \leq i_k < j_k \leq n$ such that
\begin{itemize}
\item[(i)]
$\{1,2\},\dots,\{1,n\} \in E( \mathrm{Shift}_{i_qj_q} (\mathrm{Shift}_{i_{q-1}j_{q-1}}( \cdots
(\mathrm{Shift}_{i_1j_1}(G) )\cdots )))$; \smallskip
\item[(ii)]
For $1 \leq t \leq q$, $\shift_{i_tj_t}$ is an edge shift of $\shift_{i_{t-1}j_{t-1}}(\cdots (\shift_{i_1j_1}(G)) \cdots)$.
\end{itemize}
In particular, there exist a combinatorial shifted graph $\delc G$ of $G$ such that\linebreak $m_1(\delc G)=n-1$.
\end{lemma}
\begin{proof}
We use induction on $n-\deg_G(1)$.
If $n-\deg_G(1)=1$, then $\{1,2\},\dots,\{1,n\} \in E(G)$
and $m_1 (\delc G)=n-1$ for all combinatorial shifted graph $\delc G$ of $G$ by Lemma \ref{central}.

Assume $n-\deg_G(1) >1$.
Since $G$ is connected, there exist vertices $s,t \in \nn$
such that $\{1,t\} \in E(G)$, $\{s,t\}\in E(G)$ and $\{1,s\} \not \in E(G)$.
Let $H=\shift_{1t}(G)$.

We claim that $H$ is connected.
Since $G$ is connected,
for each vertex $v \in \nn$, there is a path
$\{1,r_1\},\{r_1,r_2\},\dots,\{r_p,v\} \in E(G)$ from $1$ to $v$.
If $r_k \ne t$ for all $k=1,2,\dots,p$, then $H$ has the same path from $1$ to $v$.
If $r_k =t$ for some $k$,
then $\{1,r_{k+1}\}\in E(H)$.
Thus $\{1,r_{k+1}\},\{r_{k+1},r_{k+2}\},\dots,\{r_{p},v\}\in E(H)$ is a path in $H$ from $1$ to $v$.
Hence $H$ is connected.

Since $\{ 1,s\}\in E(H)$, we have $n-\deg_H(1)< n-\deg_G(1)$.
Thus the assertion follows from the assumption of induction.
\end{proof}

\begin{proposition}\label{chordunique}
Let $G$ be a chordal graph on $[n]$ whose combinatorial shifted graph is unique.
Then there exist integers $n_1\geq n_2\geq \dots \geq n_r>r \geq 0$ and subsets $A,B \subset [r+1,n_r]$
 with $A \cap B = \emptyset$ such that 
$G$ is isomorphic to the graph $H$ with
$$E(H)=\{\{i,j_i\}:1\leq i\leq r,\ i+1\leq j_i \leq n_i \} \cup E(K_A) \cup E(K_B).$$
\end{proposition}

\begin{proof}
\textbf{[Case 1]}
Assume $G$ has two connected components $G_1$ and $G_2$ which are not isolated vertices.
Let $A$ be the vertex set of $G_1$ and $B$ the vertex set of $G_2$.
If $E(G)=E(K_A) \cup E(K_B)$, then  $G$ has the unique combinatorial shifted graph by Lemma \ref{disj}.

Assume $E(G) \ne E(K_A) \cup E(K_B)$ and  $\{a_1,a_2\} \not \in E(G)$ for some $a_1,a_2 \in A$.
Take vertices $b_1,b_2 \in B$ with $b_1 \ne b_2$.
We may assume $a_1 < b_1$ and $a_2 < b_2$ by Lemma \ref{ape}.
Let 
$$H_1 =\shift_{a_1b_1}(G)\ \ \ \mbox{ and }  \ \ \   H_2=\shift_{a_2b_2} (H_1).$$
Then $b_1$ is an isolated vertex of $H_1$ and $H_1$ has one connected component which is not an isolated vertex.
The vertex set of this connected component is $A \cup B \setminus \{b_1\}$.
Thus, by Lemma \ref{connect},
there exists a combinatorial shifted graph $\delc {H_1}$ with 
$$m_1(\delc {H_1})=|A \cup B \setminus \{b_1\}|-1.$$
On the other hand,
for any edge $\{t,a_2\} \in E(H_1)$, we have $t \in A \setminus \{a_1\}$.
Also, for any edge $\{ t,b_2\} \in E(H_1)$,
we have $t \in B \cup \{a_1\}$.
Thus the vertex $b_2$ is an isolated vertex of $H_2$  since $A \cap B = \emptyset$.
Then $H_2$ has one connected component which is not an isolated vertex.
The vertex set of this connected component is $A \cup B \setminus \{b_1,b_2\}$.
Thus, by Lemma \ref{connect}, there exist a combinatorial shifted graph $\delc {H_2}$
such that
$$m_1 (\delc {H_2}) = |A \cup B \setminus \{b_1,b_2\}|- 1< |A \cup B \setminus \{b_1\}|- 1.$$
Since $\delc {H_1}$ and $\delc {H_2}$ are combinatorial shifted graphs of $G$,
the combinatorial shifted graph of $G$ is not unique.
Thus $G=K_A \cup K_B$.
\medskip

\textbf{[Case 2]}
If $G$ has more than three connected components which are not isolated vertices.
Then,
by applying disjoint shifts to $G$,
 we can make a graph $H$ with two connected components which are not isolated vertices.
Furthermore, we have $E(H) \ne E(K_A) \cup E(K_B)$ for all subsets $A\subset [n]$ and $ B\subset [n]$.
Thus, by [Case 1], $G$ have more than two combinatorial shifted graphs.\medskip

\textbf{[Case 3]}
Assume that $G$ has one connected component which is not an isolated vertex.
We use induction on $|E(G)|$.
If $|E(G)|=1$, then $G$ satisfies the condition of the proposition.
Assume $|E(G)|>1$.
By Lemma \ref{ape}, we may assume that there is an integer $n_1\in [n]$ such that
 $1,2,\dots,n_1$ are not isolated vertices and $n_1+1,\dots,n$ are isolated vertices of $G$.

If $E(G)=E(K_{[n_1]})$, then $G$ satisfies the condition of the proposition.
Assume $E(G) \ne E(K_{[n_1]})$.
Then, by Lemma \ref{ape}, we may assume that there exists $\{i,n_1\} \not \in E(G)$ such that  $i<n_1$.
% by Lemma \ref{ape}.
Let \begin{eqnarray*}
& &A\hspace{4pt}=\{ p \in [n_1] : \{i,p \} \in E(G),\ \{n_1,p \} \in E(G) \}, \\
& &X_i=\{ p \in [n_1] : \{i,p \} \in E(G),\ \{n_1,p \} \not \in E(G) \}, \\
& &X_{n_1}=\{ p \in [n_1] : \{i,p \} \not \in E(G),\ \{n_1,p \} \in E(G) \}, \\
& &X_0=\{ p \in [n_1] : \{i,p \} \not \in E(G),\ \{n_1,p \} \not \in E(G) \} \ \mbox{and} \ X=X_i\cup X_{n_1} \cup X_0.
\end{eqnarray*}
{\textbf [Step 1]}
We will prove $E(K_A)\subset E(G)$.
For vertices $a_1,a_2 \in A$, $G$ has the cycle 
$(\{i,a_1\},\{a_1,n_1\},\{n_1,a_2\},\{a_2,i\})$.
Since $\{i,n_1\} \not \in E(G)$ and $G$ is chordal, we have $\{a_1,a_2\} \in E(G)$.
Thus we have $E(K_A)\subset E(G)$.
\medskip \newline
{\textbf [Step 2]}
Next, we will show that there is a vertex $a \in A$ such that $a$ is a star vertex of $G$.
By [Step1], what we must prove is that 
there exists $a \in A$ such that $\{a,x\}\in E(G)$ for all $x\in X$. 
Let $A=\{a_1,a_2,\dots,a_l\}$ and $H_1=\shift_{in_1}(G)$.
Then $\{t,n_1\} \in E(H_1)$ implies $t \in A$ by the definition of $X_{n_1}$.

Assume  $a_k$ is not a star vertex of $G$ for all $k = 1,2,\dots,l$.
Then there is $x_{k} \in X$ such that $\{a_k,x_{k}\} \not \in E(G)$ for each $k=1,2,\dots,l$.
Let 
$$ H_2=\shift_{x_{1} n_1}(\shift_{x_{2} n_1}(\cdots (\shift_{x_{l} n_1}(H_1)) \cdots)).$$
Then, for each $k=1,2,\dots,l$, we have $\{a_k,n_1\} \not \in E(H_2)$  since $\shift_{x_k n_1}$ vanishes $\{a_k,n_1\}$.
Thus we have $\deg_{H_2}(n_1)=0$.
Then Lemma \ref{connect} says that there exists a combinatorial shifted graph $\delc {H_2}$ such that
$m_1(\delc{H_2})\leq n_1-2$.
However, since $G$ has the connected component on $[n_1]$,
there exists a combinatorial shifted graph $\delc G$ such that $m_1 (\delc G)=n_1-1$.
Since $G$ has the unique combinatorial shifted graph, this is a contradiction.
Thus there is a star vertex $a \in A$ of $G$.

We may assume $a=1$ by Lemma \ref{ape}.
Then Lemma \ref{central} says that
$$E(\delc G)=E(\delc {{G-\{1\}}} \cup \{\{1,2\},\dots,\{1,n_1\}\}.$$
Since $G-\{1\}$ is also a chordal graph whose combinatorial shifted graph is unique,
the assertion follows from the assumption of induction.
\end{proof}

Next, we will show some properties for edge shifts and disjoint shifts.
Let $G$ be a graph on $[n]$.
Define 
$$T_k(G)= \{A \subset [n]: |A|=k, E(K_A) \subset E(G)\}\ \ \  \mbox{ for } k =1,2,\dots,n.$$
An \textit{induced cycle} of $G$ is a cycle of $G$ which has no chords.
Thus a graph $G$ is chordal if and only if the length of every induced cycle is $3$.

\begin{lemma}\label{edgeshift}
Let $G$ be a chordal graph on $\nn$ and $\shift_{ij}$ an edge shift of $G$.
Then
\begin{itemize}
\item[(i)]  $\shift_{ij}(G)$  is also chordal;
\item[(ii)] $|T_k(G)|=|T_k({\shif})|$ for all $1 \leq k \leq n$;
\item[(iii)] If $G$ is $k$-connected, then  $\shift_{ij}(G)$ is also $k$-connected.
\end{itemize}
\end{lemma}

\begin{proof}
(i) 
Let $(\{r_0,r_1\},\{r_1,r_2\},\dots,\{r_s,r_0\})$ be an induced cycle in $\shif$.
Thus $\{r_0,r_t\}\in E(\shift_{ij}(G))$ implies $t=1$ or $t=s$.
Also, we have
\begin{eqnarray}
\{r_p,r_q\} \not \in E(G) \ \ \ \mbox{ if } \ \ \{r_p,r_q\} \subset [n] \setminus \{i,j\} \mbox{ and } 2 \leq |p-q| < s. \label{seisitu} 
\end{eqnarray}
We will show $s=2$.
If $i,j \not\in \{r_0,\dots,r_s\}$, then $(\{r_0,r_1\},\{r_1,r_2\},\dots,\{r_s,r_0\})$
is an induced cycle of $G$.
Then since $G$ is chordal, we have $s=2$.
Thus we assume $i \in \{r_0,\dots,r_s\}$ or $j \in \{r_0,\dots,r_s\}$.
\medskip

\textbf{[Case 1]}
Assume $i \in \{r_0,\dots,r_s\}$ and $i = r_0$.
If $j \in \{r_0,\dots,r_s\}$, then we have $j=r_1$ or $j= r_s$  since $\ij \in E(\shif)$.
We may assume $j=r_1$.
Then $\{j,r_2\}\in E(\shif)$.
Thus the edge $\{i,r_2\}$ must be contained in $\shif$.
Since $\{r_0,r_t\}\in E(\shif)$ implies $t=1$ or $t=s$, we have $s=2$.

If $j \not\in \{r_0,\dots,r_s\}$,
then $G$ contains one of the following cycles
\begin{eqnarray*}
&&(\{i,r_1\},\{r_1,r_2\},\dots,\{r_s,i\}),\\
&&(\{j,r_1\},\{r_1,r_2\},\dots,\{r_s,j\}),\\
&&(\{i,j\},\{j,r_1\},\{r_1,r_2\},\dots,\{r_s,i\}),\\
&&(\{j,i\},\{i,r_1\},\{r_1,r_2\},\dots,\{r_s,j\}).
\end{eqnarray*}
In each case, since $G$ is chordal,
(\ref{seisitu}) says that
$\{i,r_k\} \in E(G)$ or $\{j,r_k\} \in E(G)$ for all $k=1,2,\dots,s$.
Thus we have $\{i,r_2\}\in E(\shif)$.
Since $i=r_0$ and $\{r_0,r_t\} \in E(\shift_{ij}(G))$ implies $t=1$ or $t=s$,
we have $s=2$.
\medskip

\textbf{[Case 2]}
Assume $i \not \in \{r_0,\dots,r_s\}$, $j \in \{r_0,\dots,r_s\}$ and $j = r_0$.
Then both
$(\{i,r_1\},\{r_1,r_2\},\dots,\{r_s,i\})$ and $(\{j,r_1\},\{r_1,r_2\},\dots,\{r_s,j\})$
are cycles in $G$.
Since $G$ is chordal,
(\ref{seisitu}) says that
$\{i,r_k\} \in E(G)$ and $\{j,r_k\} \in E(G)$ for all $k=1,2,\dots,s$.
Thus we have $\{j,r_2\} \in E(\shif)$.
Since $j=r_0$, we have $s=2$ by the same way as [Case 1].
\medskip

(ii)
Fix an integer $k \in \nn$.
We will show $\shift_{ij}(T_k(G))=T_k(\shift_{ij}(G))$.
The inclusion $\shift_{ij}(T_k(G))\subset T_k(\shift_{ij}(G))$ immediately follows from
the definition of $\shift_{ij}$.
Thus what we must prove is that $ T_k(\shift_{ij}(G)) \subset \shift_{ij}(T_k(G))$.
Let $A \in T_k(\shif)$.

Assume $A \in T_k(G)$.
If $j \in A$, then, since $A \in T_k ( \shif)$,
we have
$\{a,j\} \in E(\shift_{ij}(G))$  for all $a \in A \setminus \ij$, and therefore we have
 $\{a,i\} \in E(G)$ and $\{a,j\}\in E(G)$ for all $a \in A \setminus \ij$.
This fact says that $i \in A$ or $(A \setminus \{j\}) \cup \{i\} \in T_k(G)$.
In particular, we have $C_{ij}^{T_k(G)}(A)=A \in \shift_{ij}(T_k(G))$.
On the other hand,
if $j \not\in A$, then $C_{ij}^{T_k(G)}(A)=A \in \shift_{ij}(T_k(G))$ is obvious.

Assume $A \not\in T_k(G)$.
Then, by the definition of $\shift_{ij}$, we have $i \in A$ and $j \not \in A$.
We will show $(A \setminus \{i\}) \cup \{j\} \in T_k(G)$.
Since $A \in T_k(\shif)$ and $A \not\in T_k(G)$, there is $p \in A \setminus \{i\}$ such that $\{p,i\} \not \in E(G)$ and $\{p,j \} \in E(G)$.
Also, since $A \in T_k(\shift_{ij}(G) )$,
we have $\{a,i\} \in E(G)$ or $\{a,j \}\in E(G)$  for all $a \in A \setminus \{i\}$.
For each $q \in A \setminus \{i,p\}$,
if $\{i,q\}\in E(G)$, then
the cycle $(\{i,j\},\{j,p\},\{p,q\},\{q,i\})$ is contained in $G$.
Then since $G$ is chordal and $\{p,i\} \not \in E(G)$,
we have $\{j,p\}\in E(H)$.
Thus we have $\{j , q\} \in E(G)$ for all $q \in A\setminus \{i\}$.
Hence $(A \setminus\{i\} ) \cup \{j\} \in T_k(G)$
and $C_{ij}^{T_k(G)}((A\setminus \{i\}) \cup \{j\})=A \in \shift_{ij}(T_k(G))$.

Now, we have $\shift_{ij}(T_k(G))=T_k(\shif)$.
Since $|\Gamma|= |\shift_{ij}(\Gamma)|$ holds for an arbitrary collection $\Gamma$ of subsets of $[n]$,
we have
$$|T_k(G)|=|\shift_{ij}(T_k(G))|=|T_k(\shif)|,$$
as desired.
\medskip

(iii)
Let $H=\shif$.
Assume $H$ is not $k$-connected.
Then, there exists $C=\{c_1,\dots,c_{k-1}\} \subset [n]$ such that
the induced subgraph of $H$ on $[n] \setminus C$ has connected components
$H_1, H_2, \dots,H_s$, where $s \geq 2$.
Let $A$ be the vertex set of $H_1$ and $B= [n] \setminus (A\cup C)$.
Let $L $ be the induced subgraph of $H$ on $[n] \setminus C$.

Since $G$ is $k$-connected,
there exist a vertex $a \in A$ and $b \in B$ such that $\{a,b\} \in E(G)$.
In particular, since $\{a,b\} \not \in E(H)$,
we have $a=j$ or $b=j$.
We may assume $j=a\in A$.
Since $\{j,b\} \in E(G)$, we have $\{i,b\}\in E(H)$.
Since $L$ has no edge $\{s,t\}$ with $s \in A$ and $t \in B$,
we have $i \not \in A$.
By the same way, since $\ij \in E(H)$ and $j \in A$, we have $i \not \in B$.
Hence we have $i \in C$ and $j \in A$.

We will show $\{v,i \} \not \in E(G)$ for all $v \in B$.
Let $\tilde L$ be the induced subgraph of $G$ on $[n] \setminus C$.
Since $G$ is $k$-connected,
$\tilde L$ is connected.
Since $L$ is not connected and $\tilde L$ is connected,
the vertex $j$ is the only vertex in $A$ with $\{j,t\} \in E(G)$ for some $t \in B$.
Thus, for each $v \in B$, if $\{j,r_1\},\dots,\{r_{t_v},v\}$ is a path from $j$ to $v$ 
in $\tilde L$, then we may assume $r_1,\dots,r_{t_v}\in B$.
Choose a shortest path $\{j,r_1\},\dots,\{r_{t_v},v\}$ from $j$ to $v$ for each $v \in B$.
If $\{i,v\}\in E(G)$, then $G$ has a cycle
$(\{j,r_1\},\dots,\{r_{t_v},v\},\{v,i\},\{i,j\})$.
Since $G$ is chordal
and $\{j,r_1\},\dots,\{r_{t_v},v\}$ is the shortest path from $j$ to $v$,
we have $\{i,r_k\} \in E(G)$ for all $k$.
In particular,
since $\{i,r_1\} \in E(G)$ and $\{j,r_1\} \in E(G)$,
we have $\{j,r_1\}\in E(H)$.
However, since $j \in A$ and $r_1 \in B$,
this contradicts the fact that $L$ has connected components $H_1,H_2,\dots,H_s$.

Now, we have $\{v,i\} \not \in E(G)$ for all $v \in B$.
Then, since the vertex $j\in A$ is the only vertex in $A$ with $\{j,t\} \in E(G)$ for some $t \in B$, 
we have $\{u,v\} \not \in E(G)$ for all $u\in (A \setminus\{j\}) \cup \{i\}$
and for all $v \in B$.
This fact says that the induced subgraph of  $G$ on $[n] \setminus \{(C \setminus \{i\} ) \cup \{j\}\}$
is not connected.
Since $G$ is $k$-connected, this is a contradiction.
Thus $\shif$ is $k$-connected.
\end{proof}

\begin{lemma}\label{disjoint}
Let $G$ be a graph on $[n]$ and $\shift_{ij}$ a disjoint shift of $G$.
Then
\begin{itemize}
\item[(i)]  If $G$ is chordal, then $\shif$  is also chordal.
\item[(ii)] $ |T_k(G)|=|T_k(\shif)|$ for all $k$.
\end{itemize}
\end{lemma}

\begin{proof}
Assume that $G$ has connected components $G_1,G_2,\dots,G_s$, where $s \geq 2$.
Let $A$ be the vertex set of $G_1$,
$B$ the vertex set of $G_2$.
We may assume $i \in A$ and $j\in B$.
By the definition of disjoint shifts,
we have
\begin{eqnarray}
\quad\ \ \  \{s,t\} \in E(\shift_{ij}(G)) \mbox{ and  } \{s,t\} \not \in E(G)
\Leftrightarrow
% \mbox{ if and only if }
 s=i, t\in B \mbox{ and } \{j,t\} \in E(G). \label{saig}
\end{eqnarray}
(i) Let $(\{r_0,r_1\},\dots,\{r_t,r_0\})$ be an induced cycle of $\shif$.
Since $i$ is the only vertex in $A$ which may satisfies $\{i,b\} \in E(\shif)$ of some $b \in B$,
if $i \in \{r_0,\dots,r_t\}$ then (\ref{saig}) says that
$\{r_0,\dots,r_t\} \subset A$ or $(\{r_0,\dots,r_t\} \setminus \{i\})\cup \{j\} \subset B$.
Let $i=r_0$.
If $\{i,r_1,\dots,r_t\}\subset A$, then (\ref{saig}) says that
$(\{i,r_1\},\dots,\{r_t,i\})$ is an induced cycle in $G_1$.
Since $G$ is chordal, we have $t=2$.
If $\{i,r_1,\dots,r_t\}\subset B$, then (\ref{saig}) says that
$(\{j,r_1\},\dots,\{r_t,j\})$ is an induced cycle in $G_2$, and therefore $t=2$.
In case of $i \not \in \{r_0,\dots,r_t\}$, the cycle $(\{r_0,r_1\},\dots,\{r_t,r_0\})$ is also an induced cycle in $G$.
Thus we have $t=2$.
Hence the length of every induced cycle of $G$ is $3$. 
\medskip

(ii)
Let $S \subset T_k(\shif)$.
Note that $ j\not \in S$ if $|S| \geq 2$. 
If $i \in S$, then we have $S \subset A$ or $(S \setminus \{i\}) \cup \{j\}\subset B$ by the same way as (i).
If $(S\setminus \{i\}) \cup \{j\} \subset B$, then (\ref{saig}) says that $(S \setminus \{i\}) \cup \{j\} \in T_k(G)$
and $C^{T_k(G)}_{ij}((S \setminus \{i\}) \cup \{j\})=S \in \shift_{ij}(T_k(G))$.
Also, if $i \in S$ and $S \subset A$ or if $i\not \in S$, then we have $S \subset T_k(G)$
and $C^{T_k(G)}_{ij}(S)=S \in\shift_{ij}(T_k(G))$.
These fact says that $T_k(\shif) \subset \shift_{ij} (T_k(G))$.
Since $T_k(\shif) \supset \shift_{ij} (T_k(G))$ holds for an arbitrary graph $G$,
we have $T_k(\shif) = \shift_{ij} (T_k(G))$ and $ |T_k(G)|=|T_k(\shif)|$ for all $k$.
\end{proof}

Now, we will give an algorithm.
Let $G$ be a chordal graph on $[n]$.
\medskip\newline
{\textbf{[Algorithm]}}
Let $\Delta:=\{\emptyset\}$, $V:=\nn$ and $H:=(V,E(G))$.
Repeat (I), (II) and (III) until $H=(V, \{\emptyset\})$.
\begin{itemize}
\item[(I)] 
Let $u:= \min (V)$.
If $H$ has more then 2 connected components which are not isolated vertices,
then repeat
 disjoint shifts $H:=\shift_{uj}(H)$ for some $j \in V$ with $\deg_H(j) \ne 0$ 
until $H$ has one connected components which is not an isolated vertex. 
\item[(II)]
If $H$ has one connected components which is not an isolated vertex
and $u$ is not a star vertex,
then repeat edge shifts
$H:=\mathrm{Shift}_{uv}(H)$
for all $v\in V$ with $\{u,v\} \in E(H)$
until $u$ becomes a star vertex of $H$.
\item[(III)]
If $H$ has one connected components which are not isolated vertices
and $u$ is a star vertex of $H$,
then do $\Delta := \Delta \cup \{\{u,u+1\},\dots,\{u,u+\deg_H(u)\}\}$, $V:=V \setminus \{u\}$ and $H:=(V,E(H-\{u\}))$.
\end{itemize}

Return $\Delta(G)=([n],\Delta)$.

\begin{theorem} \label{algo}
Let $G$ be a chordal graph on $[n]$.
Let $\Delta(G)$ be a graph given by the above algorithm.
Then one has $\Delta(G)=\dele G$ and  there is a combinatorial shifted graph $\delc G$ such that $\Delta(G) = \delc G$.
\end{theorem}

\begin{proof}
First, we will prove that the algorithm ends.
Since disjoint shifts reduce the numbers of connected components of $H$ which are not isolated vertices,
$H$ becomes a graph with one connected component which is not an isolated vertex by repeating (I) finitely.
Also, Lemma \ref{connect} guarantees that $u$ becomes a star vertex of $H$ by repeating (II) finitely.
Since (III) reduces $|E(H)|$, the algorithm ends.

Also, since (I) and (II) only apply $\shift_{ij}$ to $H$ for some $i,j$,
Lemma \ref{central} and (III) guarantee that
$\Delta(G)$ is a combinatorial shifted graph of $G$.

At last, we will show $\Delta(G)= \dele G$.
We claim $|T_k(G)|=|T_k(\Delta(G))|$ for all $k=1,2,\dots,n$.

We use induction on $|E(G)|$.
If $|E(G)|=1$ then $G$ is isomorphic to $\Delta(G)$ and $\dele G$.
Assume $|E(G)|>1$.
By applying (I) and (II) of the algorithm,
$G$ becomes a graph $\tilde H$ with a star vertex $1$.
Also, since (I) and (II) only use edge shifts and disjoint shifts,
Lemmas \ref{edgeshift} and \ref{disjoint} say that
$|T_k(G)|=|T_k(\tilde H)|$ for all $k$.

Consider $E( \Delta(G)-\{1\})$ and $\Delta( \tilde H-\{1\})$.
By the algorithm, we have $E( \Delta(G)-\{1\})= E(\Delta( \tilde H-\{1\}))$.
Then the induction hypothesis says that
$|T_k(\tilde H-\{1\})|=|T_k( \Delta(G) -\{1\})|$ for all $k$.
On the other hand, for every graph $G'$,
if $v$ is a star vertex of $G'$, then one has
\begin{eqnarray}
T_k(G')=T_k(G'-\{v\}) \cup \{ \{v \} \cup A: A \in T_{k-1}(G'-\{v\}) \} \mbox{ for all }k. \label{shadoho}
\end{eqnarray}
Since $|T_k(\tilde H-\{1\})|=|T_k(\Delta(G)-\{1\})|$ for all $k$ and
the vertex $1$ is a star vertex of $\tilde H$ and $\Delta(G)$,
the equality (\ref{shadoho}) says that
$|T_k(G)|=|T_k(\tilde H)|=|T_k(\Delta(G))|$ for all $k \geq 0$.

Next, we will show $I(G)$ and $I(\Delta(G))$ have the same Hilbert function.
It follows from \cite[Theorem 5.1.7]{BH} that the Hilbert function of $I(G)$ is determined by
the numbers of squarefree monomials $x_{i_1}\cdots x_{i_k}$ which do not belong to $I(G)$ for $k=1,2,\dots, n$.
Since a squarefree monomial $x_{i_1}x_{i_2}\cdots x_{i_k}$ does not belongs to $I(G)$ if and only if $\{i_1,\dots,i_k\} \in T_k(G)$,
the Hilbert function of $I(G)$ is determined by $|T_1(G)|,\dots,|T_n(G)|$.
Since we already proved $|T_k(G)|=|T_k(\Delta(G))|$ for all $k$,
it follows that $I(G)$ and $I(\Delta(G))$ have the same Hilbert function.

We recall the following easy fact:
If ideals $I$ and $J$ have a linear resolution,
then $I$ and $J$ have the same Hilbert function if and only if 
$\beta_{ij}(I)=\beta_{ij}(J)$ for all $i,j$.

Then
since both $G$ and $\Delta(G)$ are chordal,
the above fact together with Lemma \ref{codal} says that $I(G)$ and $I(\Delta(G))$ have the same graded Betti numbers.
Then, Lemma \ref{codal} also says that
$\beta_{ii+2}(I(\dele G)) = \beta_{ii+2}(I(\Delta(G)))$ for all $i$.
Since both $\Delta(G)$ and $\dele G$ are shifted,
we have $\Delta(G) = \dele G$ by Lemma \ref{kimeru}.
\end{proof}

%\begin{problem} 
%\begin{itemize} 
%\item[(i)] If $G$ is planar, then $\dele G= \dels G$?
%\item[(ii)] Prove  $m_{\leq k}(\dele G )\leq m_{\leq k}(\dels G)$ for any graph $G$ and for any $k$.
%\end{itemize}
%\end{problem}

At last, we note the following consequence about connectivity of chordal graphs.
Let $G$ be a graph on $[n]$.
We say that $G$ is \textit{$k$-hyperconnected} if $\{k,n\}\in E(\dele G)$.
In other words, $G$ is $k$-hyperconnected if and only if $\dele G$ is $k$-connected.
Also, we say that $G$ is \textit{generically $k$-rigid} if $\{k,n\} \in E(\dels G)$.
In \cite{KA-h}, Kalai proved that if $G$ is $k$-hyperconnected (generically $k$-rigid), then $G$ is $k$-connected.
In general, the converse is not true.
For example, $K_{2,2}$ is $2$-connected but $\dele {K_{2,2}}$ is not $2$-connected.
However, for chordal graphs, we have

\begin{cor}\label{connectivity}
Let $G$ be a chordal graph.
The followings are equivalent.
\begin{itemize}
\item[(i)] $G$ is $k$-connected;
\item[(ii)] $G$ is $k$-hyperconnected ( $\dele G$ is $k$-connected);
\item[(iii)] $G$ is generically $k$-rigid ( $\dels G$ is $k$-connected). 
\end{itemize}
\end{cor}

\begin{proof}
(ii) $\Leftrightarrow$ (iii) is obvious by Lemma \ref{onazi}.
(ii) $\Rightarrow$ (i) is \cite[Corollary 5.3]{KA-h}.

We will show (i) $\Rightarrow$ (ii).
We use induction on $k$.
If $k=1$, then $G$ is connected.
Thus, by Lemma \ref{connect} together with (II) of the algorithm,
we have $\{1,n\}\in \Delta (G)=\dele G$.
Assume $k>1$.
By (II) of the algorithm, we can obtain a graph $\tilde H$ with a star vertex $1$
from $G$ by applying edge shifts.
Then Lemma \ref{edgeshift} (iii) says that $\tilde H$ is also $k$-connected.
Thus $\tilde H -\{1\}$ is a $(k-1)$-connected graph on $\{2,3,\dots,n\}$.
Then the induction hypothesis says that $\{k,n\}\in E(\Delta(\tilde H - \{1\}))$.
Since $\Delta(\tilde H -\{1\})=\Delta(\tilde H)-\{1\}$, we have $\{k,n\}\in E(\Delta(\tilde H))$.
Since $\Delta(\tilde H)=\Delta (G)= \dele G$,
it follows that $G$ is $k$-hyperconnected.
\end{proof}

\end{document}